\documentclass{elsart}

\usepackage{latexsym}
\usepackage[psamsfonts]{amssymb}
\usepackage{amsfonts}
\usepackage{graphicx}
\usepackage{subfigure}
\usepackage{amsmath,cite}

\begin{document}

\newtheorem{tm}{Theorem}[section]
\newtheorem{pp}{Proposition}[section]
\newtheorem{lm}{Lemma}[section]
\newtheorem{df}{Definition}[section]
\newtheorem{tl}{Corollary}[section]
\newtheorem{re}{Remark}[section]
\newtheorem{eap}{Example}[section]

\newcommand{\pof}{\noindent {\bf Proof} }
\newcommand{\ep}{$\quad \Box$}

\newcommand{\al}{\alpha}
\newcommand{\be}{\beta}
\newcommand{\var}{\varepsilon}
\newcommand{\la}{\lambda}
\newcommand{\de}{\delta}
\newcommand{\st}{\stackrel}

\allowdisplaybreaks

\begin{frontmatter}



\title{The properties of sendograph metric on fuzzy number spaces}


\author{Huan Huang }
\
\author{}
\ead{hhuangjy@126.com (H. Huang)     }
\address
{Department of Mathematics, Jimei University, Xiamen 361021, China}

\date{}

\begin{abstract}
This paper discusses the variation of sendograph distances under some algebra operations.

\end{abstract}

\end{frontmatter}


\section{Introduction and preliminaries}

 The sendograph metric is a type of widely used metric on the fuzzy numbers spaces.
 This paper
 discusses the properties of sendograph metric.
First,
we give an overview of fuzzy number spaces,
For more details on this topic, see Ref. \cite{da,wu, huang}.

Let $\mathbb{N}$ be the set of all natural numbers, let
$\mathbb{R}^p$ be the $p$-dimensional Euclid space, and let
$F(\mathbb{R}^p)$ represent all fuzzy subsets on $\mathbb{R}^p$,
i.e. functions from $\mathbb{R}^p$ to [0, 1]. For $u\in
F(\mathbb{R}^p)$, let $[u]_{\al}$ denote the $\al$-cut of $u$, i.e.
\[
[u]_{\al}=\begin{cases}\{x\in \mathbb{R}^p : u(x)\geq \al \}, &
\al\in(0,1],    \\
{\rm supp}\, u=\overline{\{x\in \mathbb{R}^p : u(x)>0\}},&
\al=0.\end{cases}
\]
We call $ u\in F(\mathbb{R}^p)$ a $p$-dimensional fuzzy number if  $
  [u]_\al \in
K_c(\mathbb{R}^p)$ for all $ \al\in [0,1] $,
where $K_c(\mathbb{R}^p)$ denotes the set of
all nonempty compact and convex subset of $\mathbb{R}^p$. The set of all
$p$-dimensional fuzzy numbers is denoted by $E^p$.

The algebraic operations on $E^p$ are defined as follows: given $u, v\in E^p$, $\al\in [0,1]$, $r\in \mathbb{R}$,
\begin{gather*}
[u+v]_\al= [u]_\al+[v]_\al=\{x+y: x\in [u]_\al, \ y\in [v]_\al \},  \\
[r\cdot u]_\al=r\cdot [u]_\al=\{r\cdot x: x\in [u]_\al\}.
\end{gather*}

Many metrics and topologies on $E^p$ are based on the well-known
Hausdorff metric. Suppose that $K(\mathbb{R}^p)$ is the set of all nonempty
compact sets of $\mathbb{R}^p$. The {\rm Hausdorff} metric $H$ on
   $K(\mathbb{R}^p)$ is defined by:
$$H(U,V)=\max\{H^{*}(U,V),\ H^{*}(V,U)\}$$
for arbitrary $U,V\in K(\mathbb{R}^p)$, where
$$H^{*}(U,V)=\sup\limits_{u\in U}\,d\, (u,V) =\sup\limits_{u\in U}\inf\limits_{v\in
V}d\, (u,v).$$


Given $u\in E^p$, the sendograph and endograph of $u$
 are defined, respectively, by:
\begin{gather*}{\rm send}\, u =\{ (x,\al)\in [u]_0\times [0,1] : u(x)\geq
\al\},\\
\mathrm{end}\, u =\{ (x,\al)\in \mathbb{R}^p \times [0,1] : u(x)\geq
\al\}.\end{gather*} The sendograph metric $D$, the endograph metric $\Gamma$, the supremum metric $d_\infty$ and the $L_q$-type metric $d_q$ on $E^p$
are defined by
\begin{gather*}
D(u,v)=H({\rm send}\, u, {\rm send}\, v),\\
\Gamma(u,v)=H({\rm end}\, u, {\rm end}\, v),\\
d_\infty(u,v)=\sup_{\al\in [0,1]} H([u]_\al, [v]_\al),\\
d_q\, (u,v)=\left(\int_0^1 \big(d_H([u]_\al, [v]_\al) \big)^q
d\alpha  \right)^{1/q},
\end{gather*}
respectively, for all $u,v\in E^p$. It can be checked that $\Gamma(u,v)\leq D(u,v) \leq d_\infty(u,v)$ for all $u,v\in E^p$.

We say a sequence $\{u_m,\ m=1,2\ldots\}$ level converges to $u$ in $E^p$, denoted by $u_m \st{l}{\rightarrow}u$,
if
$\lim\limits_{m\to\infty}H([u_m]_\al,\ [u]_\al)=0$
 for each $\al\in[0,1]$.

Let $u, u_m$ in $E^p$, $m=1,2,\ldots$, then
$$u_m \st {d_\infty}{\rightarrow}u \Rightarrow u_m \st {l}{\rightarrow}u \Rightarrow u_m \st {D}{\rightarrow}u \Rightarrow u_m \st {d_q}{\rightarrow}u \Rightarrow u_m \st {\Gamma}{\rightarrow}u.$$

\section{The properties of sendograph metric}

In this section, we give some properties of sendograph metric on the fuzzy number spaces.

Let $x=(x_1,x_2,\ldots,x_n)\in \mathbb{R}^n$ and $\beta\in \mathbb{R}$, we use
$(x,\beta)$ to denote $(x_1,x_2,\ldots,x_n,\beta)\in \mathbb{R}^{n+1}$.

The following Theorem \ref{dso} shows that, for each fuzzy numbers $u,v,w$ in $(E^p, D)$, the distances between the points in the convex hull $\texttt{Cov}(u,v)$
 with $w$ can be dominated by the distances between $u,v$
with $w$.

\begin{tm}\label{dso}
Suppose that $u,v,w \in E^p$. Then
$$
D(\alpha u +(1-\alpha)v,w)\leq \sqrt{D(u,w)^2  +  D(v,w)^2}
$$
for each $\alpha\in [0,1]$.
\end{tm}

\pof \ Fix $\al\in [0,1]$ and put $c=\alpha u +(1-\alpha)v$.
Given $(x,\beta)\in {\rm send}\, w $, then
\begin{gather*}H((x,\beta), {\rm send}\, u)\leq D(u,w),\\
H((x,\beta), {\rm send}\, v)\leq  D(v,w)
\end{gather*}
and therefore there exists $(y_1, \lambda_1)\in {\rm send}\, u$ and
$(y_2, \lambda_2)\in {\rm send}\, v$
such that
\begin{gather*}
d((x,\beta), (y_1, \lambda_1))\leq D(u,w),\\
d((x,\beta),(y_2, \lambda_2) )\leq D(v,w) .
\end{gather*}
We may suppose that $\lambda_1 \leq \lambda_2$ with no loss of generality,
hence
$$
d((x,\beta), (\alpha y_1+ (1-\al)y_2, \lambda_1))    \leq \sqrt{D(u,w)^2 + D(v,w)^2},
$$
 and then
 $$d((x,\beta), {\rm send}\, c )\leq  \sqrt{D(u,w)^2 + D(v,w)^2}  .$$ From the arbitrariness of $(x,\beta)$ in ${\rm send}\, w$, we know
that \begin{equation}H^*({\rm send}\, w, {\rm send}\, c)\leq   \sqrt{D(u,w)^2 + D(v,w)^2}  .\label{left}\end{equation}

Given $(y,\lambda)\in {\rm send}\, c$, then there exists $(y_1,\lambda)\in {\rm send}\, u$ and $(y_2,\lambda)\in {\rm send}\, v$ such that $y=\alpha y_1+ (1-\alpha) y_2$. Since
\begin{gather*}
d((y_1,\lambda), {\rm send}\, w)\leq D(u,w), \\ d((y_2,\lambda), {\rm send}\, w)\leq D(v,w),
\end{gather*}
we can find $(x_1, \beta_1)$ and $(x_2, \beta_2)$ in ${\rm send}\, w$ such that
\begin{gather*}
d((y_1,\lambda),(x_1, \beta_1))\leq D(u,w),
\\
 d((y_2,\lambda),(x_2, \beta_2))\leq D(v,w).
 \end{gather*}
We may suppose $\beta_1 \leq \beta_2$ without loss of generality, thus
$$
d((y,\lambda),\ (\alpha x_1+ (1-\alpha) x_2, \beta_1))
\leq
\sqrt{D(u,w)^2 + D(v,w)^2},
$$
this implies that
$$
d((y,\lambda), {\rm send}\,w)
\leq
\sqrt{D(u,w)^2 + D(v,w)^2},$$
so
 we obtain that
\begin{equation}
H^*({\rm send}\, c,   {\rm send}\, w)
\leq
  \sqrt{D(u,w)^2 + D(v,w)^2} \label{r}
\end{equation}
from the arbitrariness of $(y,\lambda)$ in ${\rm send}\,c .$
Combined with inequalities \eqref{left} and \eqref{r}, we know that
$$D(c,w)\leq \sqrt{D(u,w)^2 + D(v,w)^2}. \quad \Box$$

\begin{tl} \label{dsp}
Suppose that $u,v,w \in E^p$. Then
$$
D(\alpha u +(1-\alpha)v,w)\leq \sqrt{2}\max\{D(u,w), D(v,w)\}
$$
for each $\alpha\in [0,1]$.
\end{tl}

\pof \ The desired result follows immediately from Theorem \ref{dso}. \ep

The following Theorem \ref{dou} discusses the properties of $D$-continuous fuzzy-valued-functions.

A function $F: \mathbb{R} \to E^p$ is said to be $D$-continuous if, for each $x\in \mathbb{R}$,
$\lim_{y\to x}D(F(y), F(x)) = 0$.

A set $U$ in $E^p$ is said to be uniformly-support-bounded if there is a compact set $K$ in $\mathbb{R}^p$ such that
$[u]_0 \subset K$ for all $u\in U$.
\begin{tm}\label{dou}
Let $F: \mathbb{R}\to E^p$ be $D$-continuous, then for each compact set $U$ in $\mathbb{R}$, $F(U)$ is uniformly-support-bounded.
\end{tm}

\pof \ From the basic topology, we know that $F(U)$ is a compact set in $(E^p, D)$, and thus $F(U)$ is a bounded set in $(E^p, D)$. Note that
$H([u]_0, [v]_0)\leq D(u,v)$ for all $u,v \in E^p$ , take $w\in U$, put
$$K=\{x\in \mathbb{R}^p: d(x,[F(w)]_0)\leq D(u,v)\},$$
then $K$ is a compact set in $\mathbb{R}^p$ and $K\supseteq [F(u)]_0$ for all $u\in U$,
thus $F(U)$ is uniformly-support-bounded. \ep

The following
Theorems \ref{n} and \ref{s} consider the variation of sendograph distances under the scalar multiplication operator and plus operator.

\begin{tm}\label{n}
Let $u\in E^p$ and let $\al, \beta\in \mathbb{R}$, then $D( \alpha u,  \beta u)\leq |\alpha-\beta| \max\{\|y\|: y\in [u]_0\}.$
\end{tm}

\pof \ For each $(x,\lambda)\in {\rm send}\, \alpha u$, there is a $z$ such that $x=\alpha z$
and $(z,\lambda)\in {\rm send}\, u$. Thus we know that
$$d((x,\lambda), {\rm send}\,\beta u)\leq d((\alpha z,\lambda), (\beta z,\lambda))\leq |\alpha-\beta| \max\{\|y\|: y\in [u]_0\}.$$
From the arbitrariness of $(x,\lambda)$ in ${\rm send}\, \alpha u$, we know that
$$H^*({\rm send}\, \alpha u, {\rm send}\, \beta u)\leq |\alpha-\beta| \max\{\|y\|: y\in [u]_0\},  $$
and by exchanging $\alpha$ and $\beta$ in the above inequality, we thus obtain that $$D(\alpha u, \beta u)= H({\rm send}\, \alpha u, {\rm send}\, \beta u)\leq |\alpha-\beta| \max\{\|y\|: y\in [u]_0\}. \quad \Box$$

\begin{tm} \label{s} \ Suppose that $u_j,v_j \in E^p$, $j=1,2,\ldots,n$. Let $u=\sum_{j=1}^n u_j$ and let $v=\sum_{j=1}^n v_j$, then
$D(u,v)\leq \sum_{j=1}^{n} D(u_j,v_j)$.
\end{tm}

\pof \ Given $(x,\lambda) \in {\rm send}\, u$, then there exists $x_1, x_2, \ldots, x _n$ such that $(x_j,\lambda)\in {\rm send}\, u_j$, $j=1,2,\ldots,n$ and $x=\sum_{j=1}^n x_j$. Since ${\rm send}\, v_j$ is compact, we can find
$(y_j,\mu_j)\in {\rm send}\, v_j$ such that
$d((x_j,\lambda),(y_j,\mu_j))=d((x_j,\lambda), {\rm send}\, v_j)$, $j=1,2,\ldots,n.$
Let $\mu=\min\{\mu_1,\mu_2,\ldots,\mu_n\}$,
then we have
\begin{eqnarray*}
\lefteqn{d((x,\lambda),{\rm send}\, v)\leq  d((x,\lambda),  (\sum_{j=1}^n   y_j, \mu)  )  =d((\sum_{j=1}^n x_j,\lambda),  (\sum_{j=1}^n   y_j, \mu)  ) }\\
&&\leq \sum_{j=1}^n d((x_j,\lambda), (y_j,\mu_j))= \sum_{j=1}^n d((x_j,\lambda), {\rm send}\, v_j) \\
&&\leq  \sum_{j=1}^n H^*( {\rm send} u_j,  {\rm send} v_j).
\end{eqnarray*}
From the arbitrariness of $(x,\lambda)$ in ${\rm send}\,  u$, we know that
$$H^*({\rm send}u,  {\rm send} v)\leq \sum_{j=1}^n H^*({\rm send} u_j, {\rm send} v_j). $$
Since $u$, $v$ is symmetric, we have that $$D(u,v)\leq \sum_{j=1}^n D( u_j,  v_j). \quad \Box$$

\section{Concluding remark}

This paper discusses the properties of sendograph metric on fuzzy number spaces.
The results in this paper can be used to study
the properties of fuzzy-valued functions.
It can also be used to design and analyze the fuzzy neural networks. 
Using the methods in this paper, 
we can also obtain the following properties of endograph metric.

Suppose that $u,v,w, u_j,v_j\in E^p$, $j=1,2,\ldots,n$. Then
$$\Gamma(\alpha u +(1-\alpha)v,w)\leq \sqrt{ \Gamma(u,w)^2  +  \Gamma(v,w)^2 }$$
for all $\alpha\in [0,1]$,
$$\Gamma( \beta u,   \gamma u)\leq |\beta-\gamma| \max\{\|y\|: y\in [u]_\varepsilon   \}   +\varepsilon$$
for all $\beta, \gamma \in \mathbb{R}$ and $\varepsilon \in [0,1]$, and
$$\Gamma(\sum_{j=1}^n u_j,    \sum_{j=1}^n   v_j)  \leq   \sum_{j=1}^n  \Gamma(u_j, v_j).$$



\end{document}